\documentclass[a4paper,english,leqno, 11pt]{amsart}
\usepackage{amssymb,amsfonts,amsmath,amsthm}
\usepackage{epsfig,epsf,curves}
\usepackage{bbm,epic,eepic}
\usepackage[mathscr]{eucal}
\newtheorem{thm}{Theorem}
\newtheorem{lem}{Lemma}
\newtheorem{cor}{Corollary}
\theoremstyle{remark}
\newtheorem{exmp}{Example}
\newtheorem{defn}{Definition}
\newtheorem{rem}{Remark}
\newcounter{fig}
\numberwithin{equation}{section}
\newcommand{\subbar}{\smash{\raise-7.7pt\hbox{\={}}}\mkern-8.5mu}

\newcommand{\N}{\mathbb{N}}

\newcommand{\cc}{\mathcal{C}}

\newcommand{\V}[1]{\mathit{V#1}}
\newcommand{\E}[1]{\mathit{E#1}}

\DeclareMathOperator{\Vol}{Vol}
\DeclareMathOperator{\diam}{diam}

\DeclareMathOperator{\cells}{cells}

\begin{document}

\title[Growth of self-similar graphs]{Growth of self-similar graphs}
\author[B.\ Kr\"on]{B.\ Kr\"on$^{\textstyle{\star}}$}
\begin{abstract}
Locally finite self-similar graphs with bounded geometry and without bounded geometry as well as non-locally finite self-similar graphs are characterized by the structure of their cell graphs. Geometric properties concerning the volume growth and distances in cell graphs are discussed. The length scaling factor $\nu$ and the volume scaling factor $\mu$ can be defined similarly to the corresponding parameters of continuous self-similar sets. There are different notions of growth dimensions of graphs. For a rather general class of self-similar graphs it is proved that all these dimensions coincide and that they can be calculated in the same way as the Hausdorff dimension of continuous self-similar fractals:
\[\dim X=\frac{\log \mu}{\log \nu}.\]
\end{abstract}
\bibliographystyle{plain}
\thanks{$^{\textstyle{\star}}$ The author is supported by the projects Y96-MAT and P14379-MAT of the Austrian Science Fund. Current address: Erwin Schr\"odinger Institute (ESI), Boltzmanngasse 9, 1090 - Wien, e-mail: bernhard.kroen@univie.ac.at. Mathematics Subject Classification 05C12, 28A80.}

\maketitle

\section{Introduction}

Self-similar sets are introduced in various ways. Usually they are defined as compact invariant sets of iterated function systems, confer Hutchinson \cite{hutchinson81fractals}. They are studied under different assumptions concerning their symmetries and the structure of the underlying space. Most important are the notions of nested fractals, see Lindstr{\o}m \cite{lindstroem90brownian}, and post-critically finite self-similar sets, confer Kigami \cite{kigami93harmonic}.\par
Self-similar graphs can be seen as discrete versions of these self-similar sets. There exists a lot of literature on different examples of self-similar graphs. Especially the random walk on the Sirpi\'{n}ski graph was studied extensively, see \cite{barlow88brownian}, \cite{grabner97functional1} and \cite{jones96transition}. General connections between the volume growth and the transition probabilities of the random walk were studied by Coulhon and Grigorian in \cite{coulhon98random}. Telcs studied connections between the growth dimension (also: fractal dimension), the random walk dimension and the resistance dimension in \cite{telcs89random}, \cite{telcs90spectra} and \cite{telcs95spectra}.
For a good introduction to the growth of finitely generated groups the reader is referred to the book of de la Harpe, see \cite{harpe00topics}.\par
One can define self-similarity of graphs without using a given self-similar set which is embedded into a complete metric space. A first axiomatic definition was stated by Malozemov and Teplyaev in \cite{malozemov95pure}. Their graphs correspond to fractals such that the boundaries of their \emph{cells}, see \cite{lindstroem90brownian}, contain exactly two points. With an axiomatic approach the author introduced the class of \emph{symmetrically self-similar} graphs in \cite{kroen02green}. In both papers, \cite{kroen02green} and \cite{malozemov95pure}, the spectrum of the discrete Laplacian is studied. Another approach to general self-similar graphs was chosen in \cite{malozemov01self}. In \cite{kroen02asymptotics} Teufl and the author calculated the asymptotic behaviour of the transition probabilities of the simple random walk on symmetrically self-similar graphs. They generalized results of Grabner and Woess in \cite{grabner97functional1} from the Sirpi\'{n}ski graph to these graphs.\par
Up to now, the class of symmetrically self-similar graphs is the biggest class of self-similar graphs where the simple random walk and consequently the Green functions as well as the spectrum of the Laplacian are understood well, see \cite{kroen02green} and \cite{kroen02asymptotics}. The class of graphs discussed in this note contains the class of symmetrically self-similar graphs. Several results (for example Theorems \ref{thm2_geometry} and \ref{infinite_degree} and Corollary \ref{thm_geometry2}) are relevant to these analytic studies.\par
After defining general self-similarity in Section \ref{ssgraphs} we reformulate the fixed point theorem for self-similar graphs, confer Theorem 1 in \cite{kroen02green}. This theorem can be interpreted as a graph theoretic analogue to the Banach fixed point theorem. For the more special class of homogeneously self-similar graphs, see Definition \ref{homogeneously_ss}, we discuss some basic geometric properties concerning the so-called \emph{$n$-cells}, see Definition \ref{ss_in_growth}. These $n$-cells correspond to $n$-cells and $n$-complexes in the sense of Lindstr\o m, confer \cite{lindstroem90brownian}.\par
Self-similar graphs of bounded geometry (the set of vertex degrees is bounded) correspond to finitely ramified fractals. In Section \ref{bounded} it is proved that for homogeneously self-similar graphs having a \emph{constant inner degree} (see Definition  \ref{def_inner+bounded}) there is a simple geometric equality relation between parameters, defined by the geometry of the graph, which is satisfied if and only if the graph has bounded geometry. Example \ref{exmpaustria1} shows that in general this not true for graphs without constant inner degree. The number of edges in the boundary of an \emph{$n$-cell} is calculated explicitely. We give an example of a locally finite, homogeneously self-similar graph with constant inner degree and unbounded geometry.
\par
Some basic properties of different growth dimensions are discussed in Section \ref{dimensions}.\par
In Section \ref{hgrowth} the diameter of the boundary of an $n$-cell in a homogeneously self-similar graph is computed. We give upper and lower bounds for the maximal distance between the boundary and vertices in the $n$-cell and bounds for the diameter of the whole $n$-cell. It is proved that for homogeneously self-similar graphs with bounded geometry all growth dimensions can be computed by the same formula as the Hausdorff dimension of self-similar sets which satisfy the open set condition, namely
\[\dim X=\frac{\log \mu}{\log \nu},\]
confer Hutchinson \cite{hutchinson81fractals}. Here the \emph{length scaling factor} $\nu$ is the diameter of the boundary of an $1$-cell, and the \emph{volume scaling factor} $\mu$ is the number of $1$-cells which are contained in a 2-cell. The result also holds if the diameter of a cell is greater than the length scaling factor $\nu$.
 
\section{Self-similar graphs}\label{ssgraphs}

Graphs $X=(\V X,\E X)$ with vertex set $\V X$ and edge set $\E X$ are always connected, locally finite, infinite, without loops or multiple edges. We write $\deg_X\!x$ for the \emph{degree} of a vertex $x$, which is number of vertices in $\V X$ being adjacent to $x$ in $X$. A \emph{path of length $n$} from $x$ to $y$ is an $(n+1)$-tuple of vertices
\[(z_0=x,z_1,\ldots,z_n=y)\]
such that $z_{i-1}$ is adjacent to $z_i$ for $0\le i\le n$. The distance $d_X(x,y)$ is the length of a shortest path from $x$ to $y$. A path from $x$ to $y$ is \emph{geodesic} if its length is $d_X(x,y)$. The \emph{vertex boundary} or \emph{boundary} $\theta C$ of a set $C$ of vertices in $\V X$ is the set of vertices in $\V X\backslash C$ being adjacent to some vertex in $C$. The \emph{closure} of $C$ is defined as $\overline{C}=C\cup \theta C$. Let us write $\hat C$ for the subgraph of $X$ which is spanned by the closure of $C$. We call $C$ \emph{connected} if every pair of vertices in $C$ can be connected by a path in $X$ that does not leave $C$. The set of edges $\delta C$ which connect a vertex in $C$ with a vertex in $\V X\backslash C$ is the \emph{edge boundary} of $C$.\par
For the convenience of the reader we briefly repeat the definition of self-similar graphs and their fixed point theorem, see Definitions 1 and 2 and Theorem 1 in \cite{kroen02green}.\par
Let $F$ be a set of vertices in $\V X$. Then $\cc_{X}F$ denotes the set of connected components in $\V X\setminus F$. We define the \emph{reduced graph} $X_{F}$ of $X$ by setting $\V X_{F}=F$ and connecting two vertices $x$ and $y$ in $\V X_{F}$ by an edge if and only if there exists a $C\in\cc_{X}F$ such that $x$ and $y$ are in the boundary of $C$.

\begin{defn}\label{ss_in_growth}
$X$ is \emph{self-similar} with respect to $F$ and $\psi: \V X\to \V X_F$ if
\begin{enumerate}
\item[(F1)] no vertices in $F$ are adjacent in $X$,
\item[(F2)] the intersection of the closures of two different components in $\cc_{X}F$ contains not more than one vertex and
\item[(F3)] $\psi$ is an isomorphism of $X$ and $X_F$.
\end{enumerate}
We will also write $\phi$ instead of $\psi^{-1}$, $F^{n}$ instead of $\psi^{n}F$ and we set $F^0=\V X$. Components of $\cc_{X}F^{n}$ are \emph{$n$-cells}, 1-cells are also just called \emph{cells}. The subgraphs $\hat C_n$ of $X$ which are spanned by the closures of $n$-cells are called \emph{$n$-cell graphs}, or \emph{cell graphs} instead of 1-cell graphs. An \emph{origin cell} is a cell $C$ such that $\phi\theta C\subset\overline C$. A fixed point of $\psi$ is called \emph{origin vertex}.
\end{defn}

The following lemma is a reformulation of the fixed point theorem for self-similar graphs. It is a consequence of Theorem 1 and Lemma 2 in \cite{kroen02green}.

\begin{thm}\label{banachlemma}
Let $X$ be self-similar with respect to $\tilde F$ and $\tilde \psi$. Then $X$ is also self-similar with respect to $\tilde F^k$ and $\tilde \psi^k$ for any positive integer $k$. There is an integer $n$ such that $X$, seen as self-similar graph with respect to $F=\tilde F^n$ and $\psi=\tilde\psi^n$, has either
\begin{itemize}
\item[(i)] exactly one origin cell and no origin vertex or
\item[(ii)] exactly one origin vertex $o$. And the subgraphs $X_A$ of $X$, being spanned by the closures $\overline{A}$ of components $A$ in $\cc_X\{o\}$, are self-similar graphs with respect to
\[F_A=F\cap \overline{A}\text{\quad and\quad} \psi_A=\psi|_{F_A}\]
and they have exactly one origin cell.
\end{itemize}
\end{thm}
\begin{defn}\label{homogeneously_ss}
A connected graph $X$ which is self-similar with respect to $F$ is called \emph{homogeneous} if the following axioms are satisfied:
\begin{enumerate}
\item[(H1)] All cell graphs are finite and for any pair of cells $C$ and $D$ in $\cc_{X}F$ there exists an isomorphism $\alpha: \hat C \to \hat D$ such that $\alpha \theta C=\theta D$.
\item[(H2)] Let $v_1$, $v_2$, $v_3$ and $v_4$ be vertices in the boundary $\theta C$ of a cell $C$ and $v_1\ne v_2$ and $v_3\ne v_4$, then $d_X(v_1,v_2)=d_X(v_3,v_4)$.
\end{enumerate}
   
In this section $X$ always denotes a homogeneously self-similar graph. The distance $\nu$ of two different vertices in the boundary of a cell is the \emph{length scaling factor} of $X$.
The number $\mu$ of cells in a 2-cell is called \emph{volume scaling factor} of $X$. We write $\delta_X$ instead of $|\delta C|$ and $\theta_X$ instead of $|\theta C|$ for some cell $C$ in $\cc_X F$. The diameter of a cell $C$ is denoted by $\lambda$, and we set $\rho=\lambda-\nu$.
\end{defn}

For homogeneously self-similar graphs the numbers $\lambda$, $\mu$, $\nu$, $\rho$, $\delta_X$ and $\theta_X$ are independent of the choice of the cell $C$.

\begin{exmp}\label{exmptree}
Figure \ref{baum} shows a 2-cell graph of a self-similar tree. The diameter $\lambda$ of a cell is greater than the length scaling factor $\nu$. Vertices in $F$ are drawn fat, the two vertices in $F^2$ are drawn fat and encircled.
We have $\nu=\delta_X=\theta_X=2$, $\lambda=3$ and $\mu=4$. See also Remark \ref{rem_sharp}.

\begin{center}
\begin{picture}(80,130)
\put(0,40){\circle*{2}}
\put(0,40){\circle{5}}
\put(20,40){\circle*{1}}
\put(20,60){\circle*{1}}
\put(20,80){\circle*{1}}
\put(40,40){\circle*{2}}
\put(40,80){\circle*{2}}
\put(40,100){\circle*{1}}
\put(40,120){\circle*{2}}
\put(40,60){\circle*{1}}
\put(60,60){\circle*{1}}
\put(60,100){\circle*{1}}
\put(80,100){\circle*{1}}
\put(80,60){\circle*{1}}
\put(60,40){\circle*{1}}
\put(60,20){\circle*{1}}
\put(60,0){\circle*{1}}
\put(80,40){\circle*{2}}
\put(80,40){\circle{5}}
\put(0,40){\line(1,0){80}}
\put(20,40){\line(0,1){40}}
\put(40,40){\line(0,1){80}}
\put(60,40){\line(0,-1){40}}
\put(40,60){\line(1,0){40}}
\put(40,100){\line(1,0){40}}
\refstepcounter{fig}\label{baum}
\put(30,-25){\emph{Figure \ref{baum}}}
\end{picture}
\end{center}
\end{exmp}
\vskip1cm
\begin{lem}\label{parameters}\mbox{}\par
\begin{itemize}
\item[(i)] Let $m$ and $n$ be positive integers such that $n> m$ and let $C_n$ be an $n$-cell. Then $\phi^m(C_n\cap F^{m})$ is an $(n-m)$-cell.
\item[(ii)] The number of $n$-cells in a $(n+1)$-cell $C_{n+1}$ is $\mu$ and $|\theta C_{n+1}|=\theta_X$.
\item[(iii)] Each cell graph $\hat C$ consists of $\mu$ copies of the complete graph $K_{\theta_X}$. More precisely: The image $\phi \theta C$ of the boundary of a cell $C$ spans a graph in $X$ which is isomorphic to the complete graph $K_{\theta_X}$ with $\theta_X$ vertices.
\end{itemize}
\end{lem}

\begin{proof}\mbox{}\par
\begin{itemize}
\item[(i)] The set $\theta C_n$ is the boundary of $C_n$ in $X$ as well as the boundary of $C_n\cap F^m$ in $X_{F^m}$. Since $\phi^m$ is an automorphism $X_{F^m}\to X$ the image $\phi^m \theta C_n$ is the boundary of $\phi^m(C_n\cap F^m)$ in $X$ and it is contained in $F^{n-m}$. The set $C_n \cap F^m$ is connected in $X_{F^m}$ and 
$\phi^m(C_n\cap F^m)$ is connected in $X$. It follows that $\phi^m \theta C_n$ is the boundary of the $(n-m)$-cell $\phi^m(C_n\cap F^{m})$.
\item[(ii)] For $n=1$ then the first part of the statement is clear. Suppose $n$ is greater or equal 2. Then $\phi^{n-1}(C_{n+1}\cap F^{n-1})$ is a 2-cell consisting of $\mu$ cells. These cells $C$ correspond one-to-one to the $n$-cells $D$ in $C_{n+1}$ in the following way: $\phi^n(F^n\cap D)=C$ or $F^n\cap D=\psi^n C$.\par
The image $\phi^{n+1} \theta C_{n+1}$ is the boundary of a cell, hence $|\theta C_{n+1}|=\theta_X$.
\item[(iii)] By the definition of $X_F$, the vertices in the boundary of a cell in $X$ are pairwise adjacent, thus they span a complete graph as subgraph of $X_F$. Let $C_2$ be a 2-cell in $X$. Then $\overline{C_2}\cap F$ spans $\mu$ copies of the complete graph $K_{\theta_X}$ as subgraph of $X_F$. These copies constitute a cell graph in $X_F$.
\end{itemize}
\end{proof}

\section{Bounded geometry and edge boundaries}\label{bounded}

\begin{defn}\label{def_inner+bounded}
A graph $X$ has \emph{bounded geometry} if the set of vertex degrees is bounded. A number $b$ is called \emph{constant inner degree} if $b=\deg_{\hat C}\!v$ for any vertex $v$ in the boundary of any cell $C$.
\end{defn}

\begin{thm}\label{thm1_geometry}
Let $X$ be a homogeneously self-similar graph with constant inner degree $b$, then
\[|\delta C_n|=\Big(\frac{b}{\theta_X-1}\Big)^{n-1}\delta_X\]
for any $n$-cell $C_n$.
\end{thm}

\begin{proof}
For $n=1$ the statement is clear. Let $C_n$ be an $n$-cell and let the statement of the lemma be true for $n-1$. The number of edges in $\delta(C_n\cap F)$ is $|\delta C_{n-1}|$, where $C_n\cap F$ is seen as $(n-1)$-cell in $X_F$ and $C_{n-1}$ is an arbitrary $(n-1)$-cell in $X$. Let $C$ be a cell in $X$ and let $v$ be a vertex in $F$ such that $C\subset C_n$ and $\theta (C_n\cap F)\cap \theta C=\{v\}$. Then $v$ is adjacent in $X_F$ to $\theta_X-1$ vertices in $\theta C$. Thus each cell $C$ in $C_n$ corresponds to $\theta_X-1$ edges in $\delta C_{n-1}$ and $|\delta C_{n-1}|/(\theta_X-1)$ is the number of cells $C$ in $C_n$ such that $\theta C\cap\theta C_n\ne\emptyset$. This implies
\[\delta C_n=\frac{|\delta C_{n-1}|}{\theta_X-1} b.\]
\vskip-8mm\end{proof}

\begin{thm}\label{thm2_geometry}
Let $X$ be a homogeneously self-similar graph with constant inner degree $b$. Then the following conditions are equivalent:
\begin{enumerate}
\item[(i)] $X$ has bounded geometry.
\item[(ii)] $b=\theta_X-1$.
\item[(iii)] $X$ is locally finite and $\deg_X\!v=\deg_{X_F}\!v$ for all $v\in F$.
\item[(iv)] $\delta_X=|\delta C_n|$ for any $n$-cell $C_n$.
\item[(v)] For any vertex $v$ in the boundary of any $n$-cell $C_n$ there is exactly one cell $C$ in $C_n$ such that $v\in\theta C$.
\item[(vi)] $\delta_X=\theta_X(\theta_X-1).$
\end{enumerate}
\end{thm}

\begin{proof}
The equivalence of (i), (ii) and (iii) is a slight generalization of Lemma 5 in \cite{kroen02green}, the proof stays the same. By Theorem \ref{thm1_geometry}, condition (iv) is equivalent to (ii). Condition (v) says that in any $n$-cell there are exactly $\theta_X$ different cells $C$ such that $\theta C\cap \theta C_n\ne\emptyset$. This implies $|\delta C_n|=\theta_X b$, then $X$ must have bounded geometry and $\delta_X=\theta_X(\theta_X-1)$. Condition (vi) implies $b=\theta_X-1$.
\end{proof}

As the following example shows, Theorem \ref{thm2_geometry} is in general not true for homogeneously self-similar graphs without constant inner degree.

\begin{exmp}\label{exmpaustria1}
The graph in Figure \ref{austria1} is the 4-cell graph of a homogeneously self-similar graph $X$  with bounded geometry but \[3=\delta_X>\theta_X(\theta_X-1)=2.\]
There is no constant inner degree. Vertices in $F$ are drawn fat, vertices in $F^2$ encircled, vertices in $F^3$ two times encircled and vertices in $F^4$ three times encircled.
\end{exmp}

\begin{picture}(400,200)
\setlength{\unitlength}{1.71pt}
\refstepcounter{fig}\label{austria1}
\put(0,0){\circle*{1.2}}
\put(0,0){\circle{3}}
\put(0,0){\circle{5}}
\put(0,0){\circle{7}}
\put(25,0){\circle*{1.2}}
\put(18.75,3.713084){\circle*{0.6}}
\put(12.50,0){\circle*{0.6}}
\put(37.50,10.608818){\circle*{1.2}}
\put(27.812498,9.813156){\circle*{0.6}}
\put(31.25,5.304409){\circle*{0.6}}
\put(47.187502,9.813156){\circle*{0.6}}
\put(43.75,5.304409){\circle*{0.6}}
\put(50,0){\circle*{1.2}}
\put(50,0){\circle{3}}
\put(31.25,-3.713084){\circle*{0.6}}
\put(37.50,0){\circle*{0.6}}
\put(54.781248,13.223133){\circle*{0.6}}
\put(56.25,7.577727){\circle*{0.6}}
\put(62.50,15.155455){\circle*{1.2}}
\put(51.765613,21.767020){\circle*{0.6}}
\put(59.062493,21.596523){\circle*{0.6}}
\put(75,30.310910){\circle*{1.2}}
\put(75,30.310910){\circle{3}}
\put(55.624986,28.037592){\circle*{1.2}}
\put(61.453115,33.512501){\circle*{0.6}}
\put(65.312493,29.174251){\circle*{0.6}}
\put(70.218752,17.087777){\circle*{0.6}}
\put(68.75,22.733182){\circle*{0.6}}
\put(87.50,15.155455){\circle*{1.2}}
\put(95.218752,13.223133){\circle*{0.6}}
\put(93.75,7.577727){\circle*{0.6}}
\put(98.234387,21.767020){\circle*{0.6}}
\put(90.937507,21.596523){\circle*{0.6}}
\put(94.375014,28.037592){\circle*{1.2}}
\put(88.546885,33.512501){\circle*{0.6}}
\put(84.687507,29.174251){\circle*{0.6}}
\put(79.781248,17.087777){\circle*{0.6}}
\put(81.25,22.733182){\circle*{0.6}}
\put(75,0){\circle*{1.2}}
\put(81.25,-3.713084){\circle*{0.6}}
\put(87.50,0){\circle*{0.6}}
\put(72.187502,-9.813156){\circle*{0.6}}
\put(68.75,-5.304409){\circle*{0.6}}
\put(62.50,-10.608818){\circle*{1.2}}
\put(52.812498,-9.813156){\circle*{0.6}}
\put(56.25,-5.304409){\circle*{0.6}}
\put(68.75,3.713084){\circle*{0.6}}
\put(62.50,0){\circle*{0.6}}
\put(100,0){\circle*{1.2}}
\put(100,0){\circle{3}}
\put(100,0){\circle{5}}
\put(112.50,21.650650){\circle*{1.2}}
\put(106.159374,18.094529){\circle*{0.6}}
\put(106.25,10.825325){\circle*{0.6}}
\put(105.407804,28.992924){\circle*{0.6}}
\put(111.031245,29.715517){\circle*{0.6}}
\put(109.562490,37.780384){\circle*{1.2}}
\put(115.095306,45.772182){\circle*{0.6}}
\put(117.281245,40.540842){\circle*{0.6}}
\put(125,43.301300){\circle*{1.2}}
\put(125,43.301300){\circle{3}}
\put(118.840626,25.206771){\circle*{0.6}}
\put(118.75,32.475975){\circle*{0.6}}
\put(115.939051,54.053552){\circle*{0.6}}
\put(121.562493,52.502826){\circle*{0.6}}
\put(118.124986,61.704352){\circle*{1.2}}
\put(103.531218,62.191480){\circle*{1.2}}
\put(107.032007,55.713877){\circle*{0.6}}
\put(110.828102,61.947916){\circle*{0.6}}
\put(101.703870,69.976246){\circle*{0.6}}
\put(107.390595,71.149443){\circle*{0.6}}
\put(120.310922,69.355153){\circle*{0.6}}
\put(114.687479,70.905879){\circle*{0.6}}
\put(136.157802,89.073481){\circle*{0.6}}
\put(140.312493,84.978801){\circle*{0.6}}
\put(122.906231,95.750012){\circle*{1.2}}
\put(130.266394,95.957044){\circle*{0.6}}
\put(126.765609,89.552507){\circle*{0.6}}
\put(111.249973,80.107405){\circle*{1.2}}
\put(111.249973,80.107405){\circle{3}}
\put(115.250755,93.440155){\circle*{0.6}}
\put(117.078102,87.928708){\circle*{0.6}}
\put(130.624986,83.355002){\circle*{1.2}}
\put(125.092170,77.636524){\circle*{0.6}}
\put(120.937479,81.731204){\circle*{0.6}}
\put(143.840626,68.508071){\circle*{0.6}}
\put(143.75,75.777275){\circle*{0.6}}
\put(144.592196,57.609676){\circle*{0.6}}
\put(138.968755,56.887083){\circle*{0.6}}
\put(140.437510,48.822216){\circle*{1.2}}
\put(134.904694,40.830418){\circle*{0.6}}
\put(132.718755,46.061758){\circle*{0.6}}
\put(137.50,64.951950){\circle*{1.2}}
\put(131.159374,61.395829){\circle*{0.6}}
\put(131.25,54.126625){\circle*{0.6}}
\put(200,0){\circle*{1.2}}
\put(200,0){\circle{3}}
\put(200,0){\circle{5}}
\put(200,0){\circle{7}}
\put(187.50,21.650650){\circle*{1.2}}
\put(193.840626,18.094529){\circle*{0.6}}
\put(193.75,10.825325){\circle*{0.6}}
\put(190.437510,37.780384){\circle*{1.2}}
\put(194.592196,28.992924){\circle*{0.6}}
\put(188.968755,29.715517){\circle*{0.6}}
\put(175,43.301300){\circle*{1.2}}
\put(175,43.301300){\circle{3}}
\put(184.904694,45.772182){\circle*{0.6}}
\put(182.718755,40.540842){\circle*{0.6}}
\put(181.159374,25.206771){\circle*{0.6}}
\put(181.25,32.475975){\circle*{0.6}}
\put(181.875014,61.704352){\circle*{1.2}}
\put(184.060949,54.053552){\circle*{0.6}}
\put(178.437507,52.502826){\circle*{0.6}}
\put(196.468782,62.191480){\circle*{1.2}}
\put(192.967993,55.713877){\circle*{0.6}}
\put(189.171898,61.947916){\circle*{0.6}}
\put(188.750027,80.107405){\circle*{1.2}}
\put(188.750027,80.107405){\circle{3}}
\put(198.296130,69.976246){\circle*{0.6}}
\put(192.609405,71.149443){\circle*{0.6}}
\put(179.689078,69.355153){\circle*{0.6}}
\put(185.312521,70.905879){\circle*{0.6}}
\put(150,86.602600){\circle*{1.2}}
\put(150,86.602600){\circle{3}}
\put(150,86.602600){\circle{5}}
\put(169.375014,83.355002){\circle*{1.2}}
\put(163.842198,89.073481){\circle*{0.6}}
\put(159.687507,84.978801){\circle*{0.6}}
\put(177.093769,95.750012){\circle*{1.2}}
\put(169.733606,95.957044){\circle*{0.6}}
\put(173.234391,89.552507){\circle*{0.6}}
\put(184.749245,93.440155){\circle*{0.6}}
\put(182.921898,87.928708){\circle*{0.6}}
\put(174.907830,77.636524){\circle*{0.6}}
\put(179.062521,81.731204){\circle*{0.6}}
\put(162.50,64.951950){\circle*{1.2}}
\put(156.159374,68.508071){\circle*{0.6}}
\put(156.25,75.777275){\circle*{0.6}}
\put(159.562490,48.822216){\circle*{1.2}}
\put(155.407804,57.609676){\circle*{0.6}}
\put(161.031245,56.887083){\circle*{0.6}}
\put(165.095306,40.830418){\circle*{0.6}}
\put(167.281245,46.061758){\circle*{0.6}}
\put(168.840626,61.395829){\circle*{0.6}}
\put(168.75,54.126625){\circle*{0.6}}
\put(175,0){\circle*{1.2}}
\put(181.25,-3.713084){\circle*{0.6}}
\put(187.50,0){\circle*{0.6}}
\put(162.50,-10.608818){\circle*{1.2}}
\put(172.187502,-9.813156){\circle*{0.6}}
\put(168.75,-5.304409){\circle*{0.6}}
\put(152.812498,-9.813156){\circle*{0.6}}
\put(156.25,-5.304409){\circle*{0.6}}
\put(168.75,3.713084){\circle*{0.6}}
\put(162.50,0){\circle*{0.6}}
\put(150,0){\circle*{1.2}}
\put(150,0){\circle{3}}
\put(145.218752,-13.223133){\circle*{0.6}}
\put(143.75,-7.577727){\circle*{0.6}}
\put(148.234387,-21.767020){\circle*{0.6}}
\put(140.937507,-21.596523){\circle*{0.6}}
\put(125,-30.310910){\circle*{1.2}}
\put(125,-30.310910){\circle{3}}
\put(144.375014,-28.037592){\circle*{1.2}}
\put(138.546885,-33.512501){\circle*{0.6}}
\put(134.687507,-29.174251){\circle*{0.6}}
\put(137.50,-15.155455){\circle*{1.2}}
\put(129.781248,-17.087777){\circle*{0.6}}
\put(131.25,-22.733182){\circle*{0.6}}
\put(104.781248,-13.223133){\circle*{0.6}}
\put(106.25,-7.577727){\circle*{0.6}}
\put(112.50,-15.155455){\circle*{1.2}}
\put(105.624986,-28.037592){\circle*{1.2}}
\put(101.765613,-21.767020){\circle*{0.6}}
\put(109.062493,-21.596523){\circle*{0.6}}
\put(111.453115,-33.512501){\circle*{0.6}}
\put(115.312493,-29.174251){\circle*{0.6}}
\put(120.218752,-17.087777){\circle*{0.6}}
\put(118.75,-22.733182){\circle*{0.6}}
\put(118.75,3.713084){\circle*{0.6}}
\put(112.50,0){\circle*{0.6}}
\put(125,0){\circle*{1.2}}
\put(127.812498,9.813156){\circle*{0.6}}
\put(131.25,5.304409){\circle*{0.6}}
\put(137.50,10.608818){\circle*{1.2}}
\put(147.187502,9.813156){\circle*{0.6}}
\put(143.75,5.304409){\circle*{0.6}}
\put(131.25,-3.713084){\circle*{0.6}}
\put(137.50,0){\circle*{0.6}}
\drawline(0,0)(25,0)(18.75,3.713084)(12.50,0)
(25,0)(37.50,10.608818)(27.812498,9.813156)(31.25,5.304409)
\drawline(50,0)(37.50,10.608818)(47.187502,9.813156)(43.75,5.304409)
(50,0)(25,0)(31.25,-3.713084)(37.50,0)
\drawline(50,0)(62.50,15.155455)(54.781248,13.223133)(56.25,7.577727)
(62.50,15.155455)(55.624986,28.037592)(51.765613,21.767020)(59.062493,21.596523)
\drawline(75,30.310910)(55.624986,28.037592)(61.453115,33.512501)(65.312493,29.174251)
(75,30.310910)(62.50,15.155455)(70.218752,17.087777)(68.75,22.733182)
\drawline(100,0)(87.50,15.155455)(95.218752,13.223133)(93.75,7.577727)
(87.50,15.155455)(94.375014,28.037592)(98.234387,21.767020)(90.937507,21.596523)
\drawline(75,30.310910)(94.375014,28.037592)(88.546885,33.512501)(84.687507,29.174251)
(75,30.310910)(87.50,15.155455)(79.781248,17.087777)(81.25,22.733182)
\drawline(100,0)(75,0)(81.25,-3.713084)(87.50,0)
(75,0)(62.50,-10.608818)(72.187502,-9.813156)(68.75,-5.304409)
\drawline(50,0)(62.50,-10.608818)(52.812498,-9.813156)(56.25,-5.304409)
(50,0)(75,0)(68.75,3.713084)(62.50,0)
\drawline(100,0)(112.50,21.650650)(106.159374,18.094529)(106.25,10.825325)
(112.50,21.650650)(109.562490,37.780384)(105.407804,28.992924)(111.031245,29.715517)
\drawline(125,43.301300)(109.562490,37.780384)(115.095306,45.772182)(117.281245,40.540842)
(125,43.301300)(112.50,21.650650)(118.840626,25.206771)(118.75,32.475975)
\drawline(125,43.301300)(118.124986,61.704352)(115.939051,54.053552)(121.562493,52.502826)
(118.124986,61.704352)(103.531218,62.191480)(107.032007,55.713877)(110.828102,61.947916)
\drawline(111.249973,80.107405)(103.531218,62.191480)(101.703870,69.976246)(107.390595,71.149443)
(111.249973,80.107405)(118.124986,61.704352)(120.310922,69.355153)(114.687479,70.905879)
\drawline(150,86.602600)(130.624986,83.355002)(136.157802,89.073481)(140.312493,84.978801)
(130.624986,83.355002)(122.906231,95.750012)(130.266394,95.957044)(126.765609,89.552507)
\drawline(111.249973,80.107405)(122.906231,95.750012)(115.250755,93.440155)(117.078102,87.928708)
(111.249973,80.107405)(130.624986,83.355002)(125.092170,77.636524)(120.937479,81.731204)
\drawline(150,86.602600)(137.50,64.951950)(143.840626,68.508071)(143.75,75.777275)
(137.50,64.951950)(140.437510,48.822216)(144.592196,57.609676)(138.968755,56.887083)
\drawline(125,43.301300)(140.437510,48.822216)(134.904694,40.830418)(132.718755,46.061758)
(125,43.301300)(137.50,64.951950)(131.159374,61.395829)(131.25,54.126625)
\drawline(200,0)(187.50,21.650650)(193.840626,18.094529)(193.75,10.825325)
(187.50,21.650650)(190.437510,37.780384)(194.592196,28.992924)(188.968755,29.715517)
\drawline(175,43.301300)(190.437510,37.780384)(184.904694,45.772182)(182.718755,40.540842)
(175,43.301300)(187.50,21.650650)(181.159374,25.206771)(181.25,32.475975)
\drawline(175,43.301300)(181.875014,61.704352)(184.060949,54.053552)(178.437507,52.502826)
(181.875014,61.704352)(196.468782,62.191480)(192.967993,55.713877)(189.171898,61.947916)
\drawline(188.750027,80.107405)(196.468782,62.191480)(198.296130,69.976246)(192.609405,71.149443)
(188.750027,80.107405)(181.875014,61.704352)(179.689078,69.355153)(185.312521,70.905879)
\drawline(150,86.602600)(169.375014,83.355002)(163.842198,89.073481)(159.687507,84.978801)
(169.375014,83.355002)(177.093769,95.750012)(169.733606,95.957044)(173.234391,89.552507)
\drawline(188.750027,80.107405)(177.093769,95.750012)(184.749245,93.440155)(182.921898,87.928708)
(188.750027,80.107405)(169.375014,83.355002)(174.907830,77.636524)(179.062521,81.731204)
\drawline(150,86.602600)(162.50,64.951950)(156.159374,68.508071)(156.25,75.777275)
(162.50,64.951950)(159.562490,48.822216)(155.407804,57.609676)(161.031245,56.887083)
\drawline(175,43.301300)(159.562490,48.822216)(165.095306,40.830418)(167.281245,46.061758)
(175,43.301300)(162.50,64.951950)(168.840626,61.395829)(168.75,54.126625)
\drawline(200,0)(175,0)(181.25,-3.713084)(187.50,0)
(175,0)(162.50,-10.608818)(172.187502,-9.813156)(168.75,-5.304409)
\drawline(150,0)(162.50,-10.608818)(152.812498,-9.813156)(156.25,-5.304409)
(150,0)(175,0)(168.75,3.713084)(162.50,0)
\drawline(150,0)(137.50,-15.155455)(145.218752,-13.223133)(143.75,-7.577727)
(137.50,-15.155455)(144.375014,-28.037592)(148.234387,-21.767020)(140.937507,-21.596523)
\drawline(125,-30.310910)(144.375014,-28.037592)(138.546885,-33.512501)(134.687507,-29.174251)
(125,-30.310910)(137.50,-15.155455)(129.781248,-17.087777)(131.25,-22.733182)
\drawline(100,0)(112.50,-15.155455)(104.781248,-13.223133)(106.25,-7.577727)
(112.50,-15.155455)(105.624986,-28.037592)(101.765613,-21.767020)(109.062493,-21.596523)
\drawline(125,-30.310910)(105.624986,-28.037592)(111.453115,-33.512501)(115.312493,-29.174251)
(125,-30.310910)(112.50,-15.155455)(120.218752,-17.087777)(118.75,-22.733182)
\drawline(100,0)(125,0)(118.75,3.713084)(112.50,0)
(125,0)(137.50,10.608818)(127.812498,9.813156)(131.25,5.304409)
\drawline(150,0)(137.50,10.608818)(147.187502,9.813156)(143.75,5.304409)
(150,0)(125,0)(131.25,-3.713084)(137.50,0)
\put(100,-45){\emph{Figure \ref{austria1}}}
\end{picture}
\vskip4cm
\begin{thm}\label{infinite_degree}
Let $X$ be a homogeneously self-similar graph with constant inner degree $b$ such that $b>\theta_X-1$ and let $v$ be a vertex in $\V X$. Then the following statements are equivalent:
\begin{enumerate}
\item[(i)] The degree of $v$ is infinite.
\item[(ii)] The vertex $v$ is contained in $F^n$ for any positive integer $n$.
\item[(iii)] The vertex $v$ is an origin vertex.
\end{enumerate}
\end{thm}

\begin{proof}
Let $v$ be a vertex in the boundary of an $n$-cell $C_n$. Then Theorem \ref{thm1_geometry} implies that $v$ is adjacent to
\[\frac{\delta_X}{\theta_X}\Big(\frac{b}{\theta_X-1}\Big)^{n-1}\]
vertices in $C_n$. If $v$ is in $F^n$ for any integer $n$ then it must have infinite degree. Suppose $v\in F^n\setminus F^{n+1}$. Then $\phi^n v$ is contained in $\V X\setminus F$. Since all cell graphs are finite, the number of different complete graphs $K_\theta$ which contain $v$ is finite. This is the same as the number of $n$-cells having $v$ in their boundaries. Thus $v$ has finite degree. The intersection
\[\bigcap_{n=1}^\infty F^n\]
cannot contain two different elements $x$ and $y$, because $\phi$ is a bijective contraction and $d(\phi^n x, \phi^n y)$ would tend to zero, which is impossible. Confer also Theorem \ref{nu_to_n} (i). Since $\phi F^{n+1}=F^n$ for any positive integer, we have
\[\phi\bigcap_{n=1}^\infty F^n=\bigcap_{n=1}^\infty F^n\]
and a vertex lies in this intersection if and only if it is an origin cell.
\end{proof}

As a consequence of Theorems \ref{thm2_geometry} and \ref{infinite_degree} we obtain:

\begin{cor}\label{thm_geometry2}
Let $X$ be a homogeneously self-similar graph with constant inner degree. Then one of the following statements is true:
\begin{enumerate}
\item[(i)] The graph $X$ has bounded geometry.
\item[(ii)] There exists no origin vertex and $X$ is locally finite but has unbounded geometry.
\item[(iii)] There exists an origin vertex and $X$ is non-locally finite.
\end{enumerate}
\end{cor}

\begin{exmp}\label{unbounded_graph}
The graph in Figure \ref{unbound_fig} is the 2-cell graph of a locally finite, homogeneously self-similar graph $X$ with unbounded geometry. Again, vertices in $F$ are drawn fat, vertices in $F^2$ encircled and vertices in $F^3$ two times encircled. The vertices $v_1$ and $\tilde v_1$ form the boundary of the origin cell. There is no origin vertex, $\phi v_{n+1}=v_n$ and $\phi\tilde  v_{n+1}=\tilde v_n$ for any positive integer $n$. We have $b=2$, $\theta_X=2$, $\delta_X=4$, thus $b>\theta_X-1$ and $\delta_X>\theta_X(\theta_X-1)$. Let $C_n$ be an $n$-cell and let $v_n$ be a vertex in $\theta C_n$. Then, according to Theorem \ref{thm1_geometry},
$|\delta C_n|=2^{n+1}$. And, since $v_n$ is in the boundary of three different $n$-cells, $\deg_X\! v_n=3\cdot 2^n$.
\end{exmp}
\begin{center}
\begin{picture}(85,60)
\setlength{\unitlength}{1.71pt}
\begin{drawjoin}
\jput(0,0){\circle*{1.2}}
\jput(20,7.2794){\circle*{0.6}}
\jput(20,-7.2794){\circle*{0.6}}
\jput(40,0){\circle*{1.2}}
\jput(20,7.2794){\circle*{0.6}}
\jput(20,-7.2794){\circle*{0.6}}
\jput(0,0){\circle*{1.2}}
\jput(16.304148,13.6808){\circle*{0.6}}
\jput(3.696,20.96){\circle*{0.6}}
\jput(20,34.641){\circle*{1.2}}
\jput(20,34.641){\circle{3}}
\jput(23.6959,13.6808){\circle*{0.6}}
\jput(40,0){\circle*{1.2}}
\jput(36.304,20.96){\circle*{0.6}}
\jput(23.6959,13.6808){\circle*{0.6}}
\jput(36.304,20.96){\circle*{0.6}}
\jput(20,34.641){\circle*{1.2}}
\jput(16.304148,13.6808){\circle*{0.6}}
\jput(3.696,20.96){\circle*{0.6}}
\jput(0,0){\circle*{1.2}}
\jput(16.304148,-13.6808){\circle*{0.6}}
\jput(3.696,-20.96){\circle*{0.6}}
\jput(20,-34.641){\circle*{1.2}}
\jput(23.6959,-13.6808){\circle*{0.6}}
\jput(40,-0){\circle*{1.2}}
\jput(36.304,-20.96){\circle*{0.6}}
\jput(23.6959,-13.6808){\circle*{0.6}}
\jput(36.304,-20.96){\circle*{0.6}}
\jput(20,-34.641){\circle*{1.2}}
\jput(20,-34.641){\circle{3}}
\jput(16.304148,-13.6808){\circle*{0.6}}
\jput(3.696,-20.96){\circle*{0.6}}
\jput(0,0){\circle*{1.2}}
\put(-40,0){\circle*{1.2}\circle{5}\circle{3}}
\put(-49,3){$v_3$}
\put(71,3){$\tilde v_3$}
\put(-8,3){$v_1$}
\put(42,3){$\tilde v_1$}
\put(10,35){$v_2$}
\put(10,-37){$\tilde v_2$}
\put(80,0){\circle*{1.2}\circle{5}\circle{3}}
\end{drawjoin}
\refstepcounter{fig}\label{unbound_fig}
\put(5,-50){\emph{Figure \ref{unbound_fig}}}
\end{picture}
\end{center}
\vskip3.5cm

\section{Growth dimensions}\label{dimensions}

\begin{defn}
For a vertex $x\in \V X$ and an integer $r\in\N_0$ we call
\[B(x,r)=\{y\in \V X\mid d_X(y,x)\le r\}\]
\emph{ball} (or more precisely: \emph{closed $d_X$-ball}) with centre $x$ and radius $r$. Let $A\subset \V X$ be a set of vertices. Then
\[\Vol_X\! A=\sum_{y\in A}\deg_X\!y,\]
is the \emph{volume} of $A$. We write $\Vol X$ instead of $\Vol_X\! \V X$.
\end{defn}

\begin{lem}\label{edges}
Let $X$ be any graph and let $A$ be a set of vertices in $\V X$. Then
\begin{align}
\tag{i}\Vol X&= 2|EX|\mbox{\quad and}\\
\tag{ii} \Vol\hat A &=\Vol_X\! A +|\delta A|
\end{align}
\end{lem}

\begin{proof}
In the sum of the definition of the volume each edge is counted twice.\par
In $\Vol_X\! A$ the edges connecting two vertices in $A$ are counted twice, the edges connecting a vertex in $A$ with a vertex in $\V X\setminus A$ are counted once. When we count these $|\delta A|$ edges a second time we obtain $\Vol_X\! A+|\delta A|$, the twice sum of all edges in $\E{\hat A}$, which is the same as $\Vol\hat A$.
\end{proof}
\begin{defn}
The \emph{growth function $V_x$ at $x$} is defined as
\[V_x:\quad \N_0\to\N_0\cup \{\infty\},\quad r\mapsto \Vol_X\! B(x,r).\]
We call
\[\subbar V(r)=\inf \{V_x(r)\mid x \in \V X\}\]
\emph{lower growth} or \emph{lower global growth} and
\[\bar V(r)=\sup \{V_x(r)\mid x \in \V X\}\]
\emph{upper growth} or \emph{upper global growth} of $X$.
The graph $X$ has \emph{regular volume growth}, or satisfies the \emph{doubling property}, if there exists a constant $c$ such that
\[V_x(2r)\le c\ V_x(r)\]
for any vertex $x$ and any integer $r$. We define
\[\underline\dim_G X=\liminf_{r\to\infty}\frac{\log \subbar V(r)}{\log r},\]
the \emph{lower global growth dimension,} and
\[\overline\dim_G X=\limsup_{r\to\infty}\frac{\log \bar V(r)}{\log r},\]
the \emph{upper global growth dimension} of $X$.
\end{defn}

\begin{lem}
Let $x_1$ and $x_2$ be any two vertices in a locally finite graph $Y$ of regular volume growth. Then
\[\liminf_{r\to\infty}\frac{\log V_{x_1}(r)}{\log r}=\liminf_{r\to\infty}\frac{\log V_{x_2}(r)}{\log r}\]
and
\[\limsup_{r\to\infty}\frac{\log V_{x_1}(r)}{\log r}=\limsup_{r\to\infty}\frac{\log V_{x_2}(r)}{\log r}.\]
\end{lem}

\begin{proof}
Let $r$ be an integer such that $r\ge d_X(x_1,x_2)$ and $r\ge 2$. Then
\[B(x_1,r)\subset B(x_2,d_X(x_1,x_2)+r)\subset B(x_2,2r)\]
implies
\[V_{x_1}(r)\le V_{x_2}(2r)\le c\,V_{x_2}(r)\]
and
\[\frac{\log V_{x_1}(r)}{\log r}\le\frac{\log c}{\log r}+\frac{\log V_{x_2}(r)}{\log r}.\]
\end{proof}

This lemma gives reason for the following definition:

\begin{defn}
Let $x$ be a vertex of a graph $Y$ of regular volume growth, then
\[\underline\dim X=\liminf_{r\to\infty}\frac{\log V_x(r)}{\log r}\]
is the \emph{lower growth dimension} (or \emph{lower local growth dimension}) and
\[\overline\dim X=\limsup_{r\to\infty}\frac{\log V_x(r)}{\log r}\]
is the \emph{upper growth dimension} (or \emph{upper local growth dimension}) of $X$.
\end{defn}

\begin{lem}
\[\underline\dim_G X\le\underline\dim X\le\overline\dim X\le\overline\dim_G X.\]
\end{lem}

\begin{proof}
Let $x_0$ be a vertex and $(r_n)_{n\in\N}$ be a sequence of integers such that
\[\lim_{n\to\infty}\frac{\log V_{x_0}(r_n)}{\log r_n}=\underline\dim X.\]
Then
\[\underline\dim_G X=\liminf_{r\to\infty}\frac{\log \subbar V(r)}{\log r}=
\liminf_{r\to\infty}\frac{\log \inf \{V_x(r)\mid x \in \V X\}}{\log r}\]
\[\le\liminf_{n\to\infty}\frac{\log \inf \{V_x(r_n)\mid x \in \V X\}}{\log r_n}\le
\liminf_{n\to\infty}\frac{\log V_{x_0}(r_n)}{\log r_n}=\underline\dim X.\]
The inequality relation between the upper growth dimensions follows analogously.
\end{proof}

\section{Growth of homogeneously self-similar graphs}\label{hgrowth}

In this section let $X$ always be a homogeneously self-similar graph.

\begin{thm}\label{cellvolume}
Let $C_n$ be an $n$-cell. Then
\[\Vol_X\! C_n=\Vol\hat C_n-\delta_X=\mu^n\theta_X(\theta_X-1)-\delta_X.\]
\end{thm}

\begin{proof}
By Lemma \ref{edges} (i), the volume $\Vol \hat C_n$ can be calculated by counting the edges in $\hat C_n$ twice. Let $C$ be a cell in $X$. The complete graph $K_{\theta_X}$ has $\binom{\theta_X}{2}$ edges, and Lemma \ref{parameters} (iii) implies
\[|E\hat C|=\mu\binom{\theta_X}{2}\mbox{\quad and\quad} \Vol\hat C=\mu\theta_X(\theta_X-1).\]
By Lemma \ref{parameters} (ii), $C_n$ contains $\mu$ disjoint $(n-1)$-cells $D_1,D_2,\ldots,D_\mu$ and
\[\bigcup_{k=1}^\mu \hat D_k=\hat C_n,\]
where this union means the union of graphs, not the usual set theoretic union. Thus
\[\Vol\hat C_n=\mu\Vol \hat C_{n-1}=\mu^{n-1}\Vol \hat C=\mu^n\theta_X(\theta_X-1).\]
where $C_{n-1}$ is any $(n-1)$-cell and $C$ any cell. Lemma \ref{edges} (ii) implies the rest of the statement.
\end{proof}

\begin{thm}\label{nu_to_n}
Let $C_n$ be an $n$-cell. Then
\begin{align}
\tag{i}\diam \theta C_n&=\nu^n,\\
\tag{ii}\nu^n\le\max\{d_X(x,v)\mid x\in \overline{C_n},\ v\in\theta C_n\}&\le\nu^n+\rho\frac{\nu^n-1}{\nu-1}\mbox{\qquad and}\\
\tag{iii}\nu^n\le\diam\overline{C_n}&\le\nu^n+\rho\frac{\nu^{n-1}(\nu+1)-2}{\nu-1}<\nu^{n+{\tilde\kappa}}\\\nonumber
\mbox{where\qquad} \tilde\kappa&=\frac{\log (\nu+3\rho)}{\log\nu}-1.
\end{align}
\end{thm}

\begin{proof}\mbox{}\par
\begin{itemize}
\item[(i)] By the definition of the length scaling factor, $\diam\theta C_1=\nu$. Suppose $\diam\theta C_{n-1}=\nu^{n-1}$ for all $(n-1)$-cells $C_{n-1}$.\par
Let $\pi$ be a geodesic path connecting two vertices $v$ and $w$ in the boundary $\theta C_n$. In the intersection $\pi\cap F^{n-1}$ we can find vertices $v=x_0, x_1,\ldots, w=x_n$ such that $\pi^*=(v=x_0, x_1,\ldots, w=x_n)$ is a path in $X_{F^{n-1}}$ connecting $v$ and $w$. The length of $\pi^*$ is greater or equal $\nu$. Each two consecutive vertices in $\pi^*$ are starting and end point for a path in $X$ connecting different vertices in the boundary of an $(n-1)$-cell. This means that $\pi$ decomposes into at least $\nu$ paths, each of them with length of at least $\nu^{n-1}$. Thus the length of $\pi$ is greater or equal $\nu^n$.\par
At the other hand there exists a path $\beta$ of length $\nu$ in $\overline{C_n}\cap F^{n-1}$, seen as cell in $X_{F^{n-1}}$, connecting two points in $\theta C_n$. Any pair of consecutive vertices in $\beta$ can be connected by a path in $X$ of length $\nu^{n-1}$. Thus any two points in the boundary of an $n$-cell in $X$ can be connected by a path of length less or equal $\nu^n$.
\item[(ii)] For $n=1$ we have $\nu+\rho=\lambda$. Supposed the statement is true for $n-1$. Let $\pi$ be a geodesic path connecting a vertex $v$ in $\theta C_n$ and a vertex $x$ in $\overline{C_n}$. The number of $(n-1)$-cells having vertices in common with $\pi$ is at most $\lambda$. Otherwise the\linebreak
$\phi^{n-1}$-projection of $\pi$ would be a geodesic path in a cell whose length is greater then $\lambda$. The intersection of $\pi$ with all of these $(n-1)$-cells except of the $(n-1)$-cell whose closure contains $x$ has at most length $(\lambda-1)\nu^{n-1}$. The above statement for $n-1$ says that the intersection of $\pi$ with the last cell has at most length $\nu^{n-1}+\rho\frac{\nu^{n-1}-1}{\nu-1}$. Thus the length of $\pi$ is less or equal
\[(\lambda-1)\nu^{n-1}+\nu^{n-1}+\rho\frac{\nu^{n-1}-1}{\nu-1}\]
\[=\nu^n+\rho\nu^{n-1}+\rho\frac{\nu^{n-1}-1}{\nu-1}=\nu^n+\rho\frac{\nu^n-1}{\nu-1}.\]
\item[(iii)] We can copy the proof of (ii), but we now decompose a geodesic path $\pi$ between any two vertices in $\overline{C_n}$ into at most $\lambda -2$ paths connecting two vertices in the boundary of an $(n-1)$-cell, and the initial and the end part of $\pi$. The length of the latter ones is at most $\nu^{n-1}+\rho\frac{\nu^{n-1}-1}{\nu-1}$. Thus the length of $\pi$ is less or equal
\[(\lambda-2)\nu^{n-1}+2\Big(\nu^{n-1}+\rho\frac{\nu^{n-1}-1}{\nu-1}\Big)=
\nu^n+\rho\nu^{n-1}+2\rho\frac{\nu^{n-1}-1}{\nu-1}\]
\[=\nu^n+\rho\frac{\nu^{n-1}(\nu+1)-2}{\nu-1}=\nu^n+\rho\nu^{n-1}\frac{(\nu+1)-\frac{2}{\nu^{n-1}}}{\nu-1}\]
\[<\nu^n+\rho\nu^{n-1}3=\nu^n\frac{\nu+3\rho}{\nu}.\]
Note that $\lambda=\nu+\rho$ and $\nu\ge 2$. The least real number $\tilde\kappa$ such that
\[\frac{\nu+3\rho}{\nu}\le\nu^{\tilde\kappa}\]
is
\[\tilde\kappa=\frac{\log (\nu+3\rho)}{\log\nu}-1.\]
\end{itemize}
The lower bounds in (ii) and (iii) are a consequence of (i).
\end{proof}

\begin{rem}\label{rem_sharp}
For the self-similar tree in Example \ref{exmptree} the upper bound in Lemma \ref{nu_to_n} (ii), and the first upper bound for $\diam \overline{C_n}$ in Lemma \ref{nu_to_n} (iii) are sharp.
\end{rem}

\begin{defn}
Let $\cells_X\!v$ be the number of cells $C$ such that $v$ is a vertex in $\theta C$ and let $c_X$ be
\[\sup\{\cells_X\!v\mid v\in F\}.\]
Let $M_X$ be the supremum of degrees of vertices in $\V X$. We write $c$ and $M$ instead of $c_X$ and $M_X$ if it is clear which graph is meant.
\end{defn}

The following Lemma corresponds to Lemma 4 in \cite{kroen02green}.

\begin{lem}\label{cells}
\[\cells_{X}v\ (\theta_X-1)=\deg_{X_{F}}\!v.\]
\end{lem}

\begin{cor}
\[c_X(\theta_X-1)=M_X.\]
\end{cor}

\begin{proof}
Lemma \ref{cells} implies
\[c_X(\theta_X-1)=M_{X_F}.\]
Since $X$ and $X_F$ are isomorphic $M_{X_F}$ equals $M_X$.
\end{proof}

Note that homogeneously self-similar graphs have bounded geometry if and only if $c$ is finite. Let $\kappa$ be the least integer which is greater or equal $\tilde\kappa$.

\begin{thm}\label{homogeneouslygrowth}
Let us write $r_n=\nu^n+\rho\frac{\nu^{n-1}(\nu+1)-2}{\nu-1}$ for a positive integer $n$. Then
\[r_n^{\frac{\log\mu}{\log\nu}}\theta_X(\theta_X-1)\mu^{-\kappa}\le \subbar V(r_n)\le\bar V(r_n)\]
\[\le r_n^{\frac{\log\mu}{\log\nu}}\mu^\kappa\theta_X(\theta_X-1)\big((c-1)\theta_X+1\big)+ \theta_X(\theta_X-1)(c-1)(M-1).\]
\end{thm}

\begin{proof}
According to Theorem \ref{nu_to_n} (iii) we have
\[r_n\le \nu^{n+\kappa} \mbox{\quad and\quad} n\ge \frac{\log r_n}{\log\nu}-\kappa.\]
Let $C_n$ be an $n$-cell and let $x$ be a vertex in $\overline{C_n}$. Again by Theorem \ref{nu_to_n} (iii), $C_n$ is a subset of $B(x,r_n)$. Theorem \ref{cellvolume} implies
\[\subbar V(r_n)\ge\Vol\hat C_n=\mu^n\theta_X(\theta_X-1)\ge \mu^{\frac{\log r_n}{\log\nu}-\kappa}\theta_X(\theta_X-1)=r_n^\frac{\log \mu}{\log\nu} \theta_X(\theta_X-1)\mu^{-\kappa}.\]
At the other hand let $C_{n+\kappa}$ be a $(n+\kappa)$-cell such that $x\in\overline{C_{n+\kappa}}$. Since $r_n\le \nu^{n+\kappa}$, the ball $B(x,r_n)$ is contained in the union of $C_{n+\kappa}$ and the closures of all $(n+\kappa)$-cells which are adjacent to $C_{n+\kappa}$. There are at most $(c-1)\theta_X$ of $(n+\kappa)$-cells being adjacent to $C_{n+\kappa}$. The volume of the union $D$ of $C_{n+\kappa}$ and the closures of these $(n+\kappa)$-cells is at most
\[\big((c-1)\theta_X+1\big)\mu^{n+\kappa}\theta_X(\theta_X-1)+|\delta D|,\]
the twice number of edges in the subgraph spanned by $D$, plus $|\delta D|$, confer Lemma \ref{edges} and Theorem \ref{cellvolume}. In each boundary of one of these $(n+\kappa)$-cells there are $\theta_X-1$ vertices which are not in the boundary of $C_{n+\kappa}$, and these vertices have at most $M-1$ edges in common with $\V X\setminus D$. Thus
\[|\delta D|\le (c-1)\theta_X(\theta_X-1)(M-1)\]
and
\[\bar V(r_n)\le\Vol_X D\le\big((c-1)\theta_X+1\big) \mu^{n+\kappa}\theta_X(\theta_X-1)+(c-1)\theta_X(\theta_X-1)(M-1).\]
Since $r_n\ge\nu^n$ we have
\[\mu^n\le\mu^\frac{\log r_n}{\log\nu}=r_n^\frac{\log\mu}{\log\nu}\]
and finally
\[\bar V(r_n)\le r_n^\frac{\log\mu}{\log\nu}\mu^\kappa\big((c-1)\theta_X+1\big) \theta_X(\theta_X-1)+(c-1)\theta_X(\theta_X-1)(M-1).\]
\end{proof}

The growth of a graph can be seen as the discrete analogue to the Hausdorff dimension. The main difference is that the Hausdorff dimension of sets in metric spaces depends on the underlying metric. Whereas the growth of graphs is always determined by the natural geodesic graph metric. Thus is does only depend on the subject itself.

\begin{thm}
The global lower and upper growth dimensions of homogeneously self-similar graphs of bounded geometry are
\[\underline\dim_G X=\overline\dim_G X=\frac{\log\mu}{\log\nu}.\]
\end{thm}

This means that the global growth dimensions of homogeneously self-similar graphs of bounded geometry can be obtained by the same formula as the Hausdorff dimension of self-similar sets which satisfy the open set condition, see Hutchinson \cite{hutchinson81fractals}.

\begin{proof}
For a given radius $r$ we choose an integer $n$ such that
\[\nu^n+\rho\frac{\nu^{n-1}(\nu+1)-2}{\nu-1}\le r\le\nu^{n+1}+\rho\frac{\nu^n(\nu+1)-2}{\nu-1}\]

\[=\nu\Big(\nu^n+\rho\frac{\nu^{n-1}(\nu+1)-2}{\nu-1}\Big)-\nu\rho\frac{\nu^{n-1}(\nu+1)-2}{\nu-1}+\rho\frac{\nu^n(\nu+1)-2}{\nu-1}\]
\[
=\nu\Big(\nu^n+\rho\frac{\nu^{n-1}(\nu+1)-2}{\nu-1}\Big)+\frac{\rho}{\nu-1}\big(-\nu^n(\nu+1)+2\nu+\nu^n(\nu+1)-2\big)\]
\[=\nu\Big(\nu^n+\rho\frac{\nu^{n-1}(\nu+1)-2}{\nu-1}\Big)+2\rho.\]
Then
\[\frac{r}{\nu}-\frac{2\rho}{\nu}\le\nu^n+\rho\frac{\nu^{n-1}(\nu+1)-2}{\nu-1}\le r\le\nu^{n+1}+\rho\frac{\nu^n(\nu+1)-2}{\nu-1}\le\nu r+2\rho.\]
For the radii
\[r_n=\nu^n+\rho\frac{\nu^{n-1}(\nu+1)-2}{\nu-1}\mbox{\quad and\quad}r_{n+1}=
\nu^{n+1}+\rho\frac{\nu^n(\nu+1)-2}{\nu-1}\]
we have
\[\subbar V(r_n)\le\subbar V(r)\le\bar V(r)\le\bar V(r_{n+1})\]
and by Theorem \ref{homogeneouslygrowth}
\[\Big(\frac{r}{\nu}-\frac{2\rho}{\nu}\Big)^{\frac{\log\mu}{\log\nu}}\theta_X(\theta_X-1)\mu^{-\kappa}
\le\subbar V(r)\le\bar V(r)\]
\[\le (\nu r+2\rho)^{\frac{\log\mu}{\log\nu}}\theta_X(\theta_X-1)\big((c-1)\theta_X+1\big)+ (c-1)\theta_X(\theta_X-1)(M-1)\]
for any integer $r$. It follows that
\[\liminf_{r\to\infty}\frac{\log \subbar V(r)}{\log r}=\limsup_{r\to\infty}\frac{\log \bar V(r)}{\log r}=\frac{\log\mu}{\log\nu}.\]
\end{proof}\par\noindent
\emph{Remark.} This paper is based on parts of the author's PhD thesis \cite{kroen01spectral}.\\ \par\noindent
\emph{Acknowledgement.} The coordinates for the `Austria'-graph (fractal mountains looking like the shape of the country on a map) in Figure \ref{austria1} where computed by a program for visualizing self-similar graphs which was written by E.~Teufl.


\begin{thebibliography}{10}

\bibitem{barlow88brownian}
M.~T. Barlow and E.~A. Perkins.
\newblock Brownian motion on the {S}ierpi\'nski gasket.
\newblock {\em Probab. Theory Related Fields}, 79(4):543--623, 1988.

\bibitem{coulhon98random}
T.~Coulhon and A.~Grigoryan.
\newblock Random walks on graphs with regular volume growth.
\newblock {\em Geom. Funct. Anal.}, 8(4):656--701, 1998.

\bibitem{harpe00topics}
P.~de~la Harpe.
\newblock {\em Topics in geometric group theory}.
\newblock University of Chicago Press, Chicago, IL, 2000.

\bibitem{grabner97functional1}
P.~J. Grabner and W.~Woess.
\newblock Functional iterations and periodic oscillations for simple random
  walk on the {S}ierpi\'nski graph.
\newblock {\em Stochastic Process. Appl.}, 69(1):127--138, 1997.

\bibitem{hutchinson81fractals}
J.~E. Hutchinson.
\newblock Fractals and self-similarity.
\newblock {\em Indiana Univ. Math. J.}, 30(5):713--747, 1981.

\bibitem{jones96transition}
O.~D. Jones.
\newblock Transition probabilities for the simple random walk on the
  {S}ierpi\'nski graph.
\newblock {\em Stochastic Process. Appl.}, 61(1):45--69, 1996.

\bibitem{kigami93harmonic}
J.~Kigami.
\newblock Harmonic calculus on p.c.f.\ self-similar sets.
\newblock {\em Trans. Amer. Math. Soc.}, 335(2):721--755, 1993.

\bibitem{kroen01spectral}
B.~Kr{\"o}n.
\newblock {\em Spectral and structural theory of infinite graphs}.
\newblock PhD thesis, Graz University of Technology, 2001.

\bibitem{kroen02green}
B.~Kr{\"o}n.
\newblock Green functions on self-similar graphs and bounds for the spectrum of
  the {L}aplacian.
\newblock to appear in \emph{Ann.~Inst.~Fourier} \textbf{52}, no. 6, 2002.

\bibitem{kroen02asymptotics}
B.~Kr{\"o}n and E.~Teufl.
\newblock Asymptotics of the transition probabilities of the simple random walk
  on self-similar graphs.
\newblock preprint, 2002.

\bibitem{lindstroem90brownian}
T.~Lindstr{\o}m.
\newblock Brownian motion on nested fractals.
\newblock {\em Mem. Amer. Math. Soc.}, 83(420):iv+128, 1990.

\bibitem{malozemov95pure}
L.~Malozemov and A.~Teplyaev.
\newblock Pure point spectrum of the {L}aplacians on fractal graphs.
\newblock {\em J. Funct. Anal.}, 129(2):390--405, 1995.

\bibitem{malozemov01self}
L.~Malozemov and A.~Teplyaev.
\newblock Self-similarity, operators and dynamics.
\newblock preprint, 2001.

\bibitem{telcs89random}
A.~Telcs.
\newblock Random walks on graphs, electric networks and fractals.
\newblock {\em Probab. Theory Related Fields}, 82(3):435--449, 1989.

\bibitem{telcs90spectra}
A.~Telcs.
\newblock Spectra of graphs and fractal dimensions. {I}.
\newblock {\em Probab. Theory Related Fields}, 85(4):489--497, 1990.

\bibitem{telcs95spectra}
A.~Telcs.
\newblock Spectra of graphs and fractal dimensions. {I}{I}.
\newblock {\em J. Theoret. Probab.}, 8(1):77--96, 1995.

\end{thebibliography}
\end{document}